
%
\input amstex
\documentstyle{amsppt}
\magnification=\magstephalf
\NoRunningHeads
\refstyle{A}
\widestnumber\key{SSSSSSS}   
 \addto\tenpoint{\baselineskip 15pt
  \abovedisplayskip18pt plus4.5pt minus9pt
  \belowdisplayskip\abovedisplayskip
  \abovedisplayshortskip0pt plus4.5pt
  \belowdisplayshortskip10.5pt plus4.5pt minus6pt}\tenpoint
\pagewidth{6.5truein} \pageheight{8.9truein}
\subheadskip\bigskipamount
\belowheadskip\bigskipamount
\aboveheadskip=3\bigskipamount
\catcode`\@=11
\def\output@{\shipout\vbox{%
 \ifrunheads@ \makeheadline \pagebody
       \else \pagebody \fi \makefootline
 }%
 \advancepageno \ifnum\outputpenalty>-\@MM\else\dosupereject\fi}
\outer\def\subhead#1\endsubhead{\par\penaltyandskip@{-100}\subheadskip
  \noindent{\subheadfont@\ignorespaces#1\unskip\endgraf}\removelastskip
  \nobreak\medskip\noindent}
\outer\def\enddocument{\par
  \add@missing\endRefs
  \add@missing\endroster \add@missing\endproclaim
  \add@missing\enddefinition
  \add@missing\enddemo \add@missing\endremark \add@missing\endexample
 \ifmonograph@ 
 \else
 \vfill
 \nobreak
 \thetranslator@
 \count@\z@ \loop\ifnum\count@<\addresscount@\advance\count@\@ne
 \csname address\number\count@\endcsname
 \csname email\number\count@\endcsname
 \repeat
\fi
 \supereject\end}
\catcode`\@=\active
\CenteredTagsOnSplits
\NoBlackBoxes
\def\today{\ifcase\month\or
 January\or February\or March\or April\or May\or June\or
 July\or August\or September\or October\or November\or December\fi
 \space\number\day, \number\year}
\define\({\left(}
\define\){\right)}
\define\Ahat{{\hat A}}
\define\Aut{\operatorname{Aut}}
\define\CC{{\Bbb C}}

\define\RR{{\Bbb R}}

\define\ZZ{{\Bbb Z}}
\define\[{\left[}
\define\]{\right]}

\define\chiup{\raise.5ex\hbox{$\chi$}}

\define\exertag #1#2{\removelastskip\bigskip\medskip\eightpoint\noindent%
\hbox{\rm\ignorespaces#2\unskip} #1.\ }

\define\inv{^{-1}}

\define\protag#1 #2{#2\ #1}

\define\res#1{\negmedspace\bigm|_{#1}}
\define\temsquare{\raise3.5pt\hbox{\boxed{ }}}

\define\theprotag#1 #2{#2~#1}

\define\zmod#1{\ZZ/#1\ZZ}

\def\bigstrut{\hbox{\vrule height18pt depth 9pt width0pt}}
\def\entry#1:#2:#3 {\bigstrut{#1}&{#2}&{#3}\cr}
\def\titlestrut{\hbox{\vrule height12pt depth 9pt width0pt}}
\define\RZ{\RR/\ZZ}
\define\Sect{\operatorname{Sect}}
\define\TT{\Bbb T}
\define\Trace{\operatorname{Trace}}
\define\Zero{\operatorname{Zero}}
\define\bK{\partial K}
\define\bL{\partial L}
\define\bP{\partial P}
\define\bQ{\partial Q}
\define\bX{\partial X}
\define\bY{\partial Y}
\define\conn#1{\Cal A_{#1}}
\define\cut{^{\text{cut}}}
\define\el#1{e^{2\pi i\lambda _X(#1)}}
\define\fld#1{\Cal C_{#1}}
\define\fldb#1{\overline{\Cal C_{#1}}}
\define\hol{\operatorname{hol}}
\define\sE{\Cal E}
\define\sG{\Cal G}
\define\sM{\Cal M}
\define\sgn{\operatorname{sgn}}
\define\tXL{\tors{\bX}{\bL}}
\define\tors#1#2{T_{#1}(#2)}
\input epsf           
	\topmatter
 \title\nofrills Characteristic Numbers and Generalized Path Integrals
\endtitle
 \author Daniel S. Freed  \endauthor
 \thanks The author is supported by NSF grant DMS-9307446, a Presidential
Young Investigators award DMS-9057144, and by the O'Donnell Foundation.  He
thanks the University of Warwick and the Geometry Center at the University of
Minneapolis for their hospitality while this note was written.\endthanks
 \affil Department of Mathematics \\ University of Texas at Austin\endaffil
 \address Department of Mathematics, University of Texas, Austin, TX
78712\endaddress
 \email dafr\@math.utexas.edu \endemail
 \date November 30, 1993\enddate
 \dedicatory To Raoul Bott in honor of his $70^{\text{th} }$ birthday
\enddedicatory
	\endtopmatter

\document

I hope in this exposition to touch on two themes which are among Professor
Bott's many mathematical interests.  The first---characteristic
numbers---appears throughout his work.  My second topic, path integrals, have
been around in physics since they were introduced by Feynman in the late
1940s, but have only recently been applied to problems of purely geometric
interest.  Edward Witten has led the way in this ``topological quantum field
theory'', which has attracted the enthusiasm of many mathematicians.  The
ideas of quantum field theory, though far from mathematically understood,
seem to provide a unified framework for many recent invariants in low
dimensional topology (\`a la Donaldson, Jones, Casson, Floer, \dots) and
introduce some new invariants as well.  Our first goal is to explain that
ideas of (classical) field theory inspire new insight into characteristic
numbers.  The second goal is to explain how summing over fields---the path
integral---produces diffeomorphism invariants which satisfy {\it gluing
laws\/}, and how a generalization of this idea explains some of the algebraic
structure behind these invariants in three dimensions.  We can only do this
rigorously in a ``toy model'', where the path integral reduces to a finite
sum, but we hope that the ideas here will shed light on more interesting
examples as well.  The toy model is based on a characteristic number for
finite principal covering spaces.  As I have explained most of these ideas
elsewhere, my goal here is to give an elementary account of the basic
concepts and the simplest cases.

The oldest, and in some sense most intriguing, characteristic number is the
Euler number of a manifold.  More generally, an oriented real vector bundle
of rank~$n$ over a closed oriented manifold of dimension~$n$ has an Euler
number, which takes values in the integers.  Real and complex vector bundles
over closed oriented manifolds have other characteristic numbers due to
Chern, Pontrjagin, and Stiefel-Whitney.  These integers are {\it primary\/}
characteristic numbers in {\it cohomology\/}; they are topological
invariants.  There are similar {\it primary\/} characteristic numbers in {\it
$K$-theory\/} for real and complex bundles over compact oriented (or spin)
manifolds.  The topological definition~\cite{AH} of these characteristic
numbers uses Bott Periodicity in an essential way.  They include the
signature of a manifold, the $\Ahat$-genus, etc.  What I term {\it secondary
characteristic numbers\/} are {\it geometric\/} invariants which depend on a
metric or connection and take values in~$\RZ$, or, after exponentiating, in
the circle group~$\TT$.  In cohomology they are the Chern-Simons numbers and
in $K$-theory they are the $\eta $-invariants.  We summarize in Table~1.  Our
observation is that these classical invariants extend to invariants of
manifolds ``below the top dimension'' and of manifolds with boundary.  Thus,
for example, the characteristic number~$c_4$ is an integer associated to a
complex vector bundle over a closed oriented 8-manifold.  There are
invariants of complex vector bundles over compact oriented manifolds
(possibly with boundary) in all dimensions~$\le 8$, but now the invariants
are not simply numbers.  For manifolds of dimension~7 they are sets, for
manifolds of dimension~6 they are categories, and so on.  These sets and
categories have more structure, of course.  These extended invariants obey
gluing laws, so are in some sense {\it local\/} invariants.  That is, in some
respect they behave in some respects like the integral over the manifold of a
locally computed density.  In~\S{1} we explain the simplest case of these
ideas: the Euler number of a complex line bundle over a surface.  The
structure of the invariants in other cases is similar, though the actual
constructions are quite different.  I do not know all of the constructions
for the $K$-theory invariants.

\bigskip
\midinsert
$$ \vcenter{  \offinterlineskip \tabskip = 0pt \halign{
	\vrule\enspace\hfil#\hfil\enspace\vrule\hskip2pt
        \vrule&\enspace\hfil#\hfil\enspace
	&\vrule\enspace\hfil#\hfil\enspace\vrule\cr
\noalign{\hrule}
\titlestrut&{\bf cohomology}&{\bf $K$-theory} \cr
\noalign{\hrule}
\vphantom{\vrule height 2pt}&&\cr \noalign{\hrule}
\entry {\bf primary}:{Chern numbers, Pontrjagin numbers}:{signature,
$\Ahat$-genus}
\noalign{\hrule}
\entry {\bf secondary}:{Chern-Simons invariants}:{$\eta $-invariants}
\noalign{\hrule}
}} $$
\nobreak
 \centerline{Table~1: Classical characteristic numbers}
 \medskip
\endinsert

Because the invariants in Table~1 obey gluing laws, they are appropriate
classical actions for a field theory.  Whereas the classical action is an
``integral'' over a finite dimensional manifold, the quantum path integral is
an ``integral'' over a space of fields on that manifold, which typically is
infinite dimensional.  Although mathematical physicists have made great
progress in understanding the technology of these integrals, the cases of
topological significance discussed here remain elusive.  Our main idea
in~\S{2} is to extend this (ill-defined) notion of integration over fields to
the {\it extended\/} classical action, i.e., to fields on manifolds below the
top dimension.  For example, in a three dimensional topological field theory
the path integral for a closed oriented 3-manifold produces a complex number,
whereas canonical quantization assigns a complex Hilbert space to a closed
oriented 2-manifold.  These are usually considered as different processes.
We reinterpret this Hilbert space as the result of an integration process
over the space of fields on the 2-manifold.  The extended classical action
has values which {\it are\/} hermitian lines, and the integral (direct sum)
of hermitian lines is a Hilbert space.  This is very strange---integrals are
usually numbers, not Hilbert spaces!  Nonetheless, this strange notion of
integration leads immediately to a gluing law for the Hilbert spaces which,
for example, is closely related to the ``Verlinde formula'' in conformal
field theory.  This idea is a pure formality in general; we emphasize that so
far this whole scheme only works rigorously for our toy model.  After all,
the usual path integral is also not rigorously defined in topological
theories.  Our extended notion of path integral applies to 1-manifolds (in a
three dimensional topological theory) and leads to a braided monoidal
category, and so ultimately to a quantum group.  We view this as the {\it
solution\/} to the theory.  That is, we start with the classical action to
define the theory, then compute the quantum group (actually its category of
representations) from the theory, and finally derive formulas for invariants
from the gluing laws of the path integral (e.g.~in terms of a presentation as
surgery on a link).

\bigskip \midinsert
$$ \vcenter{  \offinterlineskip \tabskip = 0pt \halign{
	\vrule\enspace\hfil#\hfil\enspace\vrule\hskip2pt
        \vrule&\enspace\hfil#\hfil\enspace
	&\vrule\enspace\hfil#\hfil\enspace\vrule\cr
\noalign{\hrule}
\titlestrut&{\bf cohomology}&{\bf $K$-theory} \cr
\noalign{\hrule}
\vphantom{\vrule height 2pt}&&\cr \noalign{\hrule}
\entry {\bf primary}:{Donaldson invariants, Floer homology, Gromov
invariants}:{???}
\noalign{\hrule}
\entry {\bf secondary}:{Reshetikhin-Turaev-Witten invariants, Jones
invariants}:{???}
\noalign{\hrule}
}} $$
\nobreak
 \centerline{Table~2: Quantum characteristic numbers}
 \medskip
\endinsert

I find it useful to organize many of the topological quantum field theories
floating around into Table~2.  I like to think of the invariants in these
theories as {\it quantum characteristic numbers\/}.  Often people distinguish
between two types of topological field theories, which here is the
distinction between primary and secondary invariants.  Notice that on the
quantum level the geometric data used to define the secondary characteristic
number has been integrated out, so we are left with topological invariants in
both types of theory.  The main example here is Witten's definition~\cite{W1}
of an invariant of 3-manifolds as a path integral of Chern-Simons.  This is
most closely analogous to the usual path integral in physics.  Less clear is
the integration process which leads from primary characteristic numbers
(which are already topological invariants) to Donaldson invariants, etc.
Witten introduced {\it supersymmetric\/} path integrals~\cite{W2}, ~\cite{W3}
to explain these invariants, and Baulieu-Singer reinterpreted this in terms
of primary characteristic numbers and gauge fixing~\cite{BS}.  We should also
mention that these invariants are in some sense an infinite dimensional Euler
number, and that the path integral representation is an infinite dimensional
version of the Mathai-Quillen formula~\cite{MQ}, ~\cite{AJ}.  From this point
of view the primary characteristic numbers in $K$-theory which I placed in
Table~1 could be pigeon-holed in the upper left hand corner of Table~2,
thanks to supersymmetric quantum mechanics.  Just as with the classical
invariants in Table~1, the geometric and topological applications of {\it
primary\/} quantum invariants (Donaldson) have been more striking than the
applications of {\it secondary\/} invariants (Chern-Simons).

The idea that in a three dimensional topological quantum field theory one
should attach certain types of categories to 1-manifolds has been discussed
by Kazhdan, Segal, Lawrence, Yetter, Crane, and many others.  The
construction by generalized path integrals is new.  I refer the reader
to~\cite{F1} for more details about generalized path integrals and for the
detailed computations for finite group gauge theory.  The appendix
to~\cite{FQ} defines the characteristic numbers used in this model, and
\cite{FQ}~contains many more details about the basics of the theory.  There
is another expository account of some of this material in~\cite{F2}.  In
particular, the derivation of the quantum group is explained in more detail
there.

It is a great pleasure and honor to dedicate this paper to Raoul Bott.  His
birthday is a wonderful occasion to celebrate his youthful exuberance for
life and for mathematics.  {\it L'chaim\/}, Raoul!

\newpage
\head
\S{1} Characteristic Numbers
\endhead
\comment
lasteqno 1@ 13
\endcomment

Consider a complex line bundle $L\to X^2$ over a closed oriented
2-manifold~$X$.  Its Euler number~$e_X(L)\in \ZZ$ can be computed in several
ways.  For example, choose any section $s\:X\to L$ which is transverse to the
zero section.  Then $s$~has isolated zeros, and the orientations determine a
sign~$\sgn(y)$ for each zero~$y$.  A standard argument shows that
  $$ e_X(L) = \sum\limits_{y\in \Zero(s)}\sgn(y) \tag{1.1} $$
is independent of~$s$, and so defines a topological invariant $e_X(L)\in
\ZZ$.  As we deform~$s$ two zeros of opposite sign can simultaneously die or
simultaneously appear.  But the total signed number of zeros is constant.

Consider now the same situation over a compact oriented 2-manifold~$X$ with
nonempty boundary.  Assume that the section~$s$ does not vanish on~$\bX$.  It
is still true that \thetag{1.1}~does not change if we modify~$s$ in the
interior of~$X$, but now \thetag{1.1}~depends on the restriction of~$s$ to
the boundary.  If $s$~is deformed allowing zeros on the boundary, then
\thetag{1.1}~changes according to the number of zeros (counted with sign, of
course) which flow through the boundary during the deformation.  More
formally, define $e_X(L,t)$ for any nonvanishing section $t\:\bX\to L$ by
extending to a section $s\:X\to L$ and counting zeros as in~\thetag{1.1}.
Then for two such section~$t,t'$ we have
  $$ e_X(L,t') - e_X(L,t) = \deg_{\bX}(t'/t). \tag{1.2} $$
Here the ratio~$t'/t$ is a map $\bX\to\CC^\times $ which has a degree, or
winding number, computed using the induced orientation of~$\bX$.  So the
relative Euler number depends on the section~$t$ (up to homotopy).  This is
the traditional point of view on relative characteristic numbers---they
depend on a trivialization on the boundary.

Now the new twist: We can define a relative Euler number which is independent
of~$t$ if we abandon the idea that it should be a number.  Namely, let
$\Sect^*$~denote the set of nonzero sections, and set
  $$ \tXL = \{e\:\Sect^*(\partial L)\longrightarrow \ZZ \text{ which obey
     $e(t')-e(t)= \deg_{\bX}(t'/t)$} \}. \tag{1.3} $$
Notice that if $e,e'\in \tXL$ then $e'-e$~is a constant integer.  Likewise,
if $e\in \tXL$ and $n\in \ZZ$ then $e+n\in \tXL$.  So $\tXL$~is a principal
homogeneous space for the integers, a so-called {\it $\ZZ$-torsor\/}.  It is
an ``affine'' copy of~$\ZZ$---a copy of~$\ZZ$ without a preferred origin.  In
other terms it is a principal $\ZZ$~bundle over a point.  It is {\it not\/} a
group.  The function~$e_X(L,\cdot )$ lies in this torsor, by~\thetag{1.2},
and so we obtain a relative Euler number
  $$ e_X(L)\in \tXL \tag{1.4} $$
which only depends on the bundle~$L$.

We now explain the sense in which \thetag{1.4}~is a topological invariant.
First, notice that~\thetag{1.3} defines a $\ZZ$-torsor~$\tors YK$ for any
line bundle $K\to Y$ over any closed oriented 1-manifold~$Y$.  It is a
topological invariant of~$K$.  Usually topological invariants are {\it
numbers\/}, and equivalent objects have {\it equal\/} invariants.  Here the
topological invariant is a {\it set\/}, and equivalent objects have {\it
isomorphic\/} invariants.  Precisely, if $K' @>\psi >> K$ is an isomorphism
of line bundles over~$Y$, it induces an isomorphism of $\ZZ$-torsors $\tors
Y{K'} @>\psi _*>> \tors YK$.  This is familiar from algebraic topology, where
invariants like the fundamental group or homology groups of a space are sets
and homeomorphisms on spaces induce isomorphisms of these sets.  The relative
Euler number~\thetag{1.4} is a topological invariant in this sense: If $L'
@>\varphi >> L$ is an isomorphism of circle bundles over~$X$, then $(\partial
\varphi _*)\bigl(e(L') \bigr) = e(L)$.  Most importantly, the Euler number
obeys a {\it gluing law\/}.  To see this, first note that if $-Y$~denotes
$Y$~with the opposite orientation, then for any $K\to Y$ there is a natural
pairing
  $$ +\:\tors YK\times \tors{-Y}K\longrightarrow \ZZ \tag{1.5} $$
by addition.  Now suppose $X_1,X_2$~are oriented surfaces with boundaries
$\partial X_1=Y$ and $\partial X_2=-Y$.  Let $X=X_1\cup_Y X_2$ denote the
surface formed by gluing along the boundary.  Suppose we are also given
circle bundles~$L_i\to X_i$ and an isomorphism $\partial L_1\to \partial
L_2$.  Let $L\to X$ denote the glued bundle.  Then the gluing law asserts
  $$ e_X(L) = e_{X_1}(L_1) + e_{X_2}(L_2). \tag{1.6} $$
The proof is direct from the definitions.

There are alternative constructions of these invariants.  For example, fix a
metric on~$L$ and let $\conn L$~denote the space of unitary connections
on~$L$.  Let $F(\theta )\in \Omega ^2_X$ denote $i/2\pi $~times the curvature
of a connection~$\theta \in \conn L$.  Then
  $$ e_X(L,\theta ) = \int_{X}F(\theta ) \tag{1.7} $$
is a real-valued function on~$\conn L$.  If $\bX=\emptyset $ it is a constant
integer equal to~$e_X(L)$.  If $\bX\not= \emptyset $ then
  $$ \exp\bigl(2\pi i\int_{X}F(\theta ) \bigr) = \hol_{\bX}(\theta ),
      $$
where $\hol_{\bX}(\theta )$~is the holonomy of~$\theta $ around~$\bX$.  In this
situation we define the $\ZZ$-torsor
  $$ T'_{\bX}(\bL) = \{f\: \conn{\partial L}\to\RR : e^{2\pi if(\theta )} =
     \hol_{\bX}(\theta )\}. \tag{1.8} $$
Then \thetag{1.7}~defines the relative Euler number.  With these definitions
the gluing law~\thetag{1.6} follows from the additivity of integration:
  $$ \int_{X} = \int_{X_1} + \int_{X_2}.  $$

The Euler number is an example of a primary characteristic number (see
Table~1).  More generally, if $G$~is a Lie group and $\lambda \in
H^n(BG;\ZZ)$ a universal characteristic class, then $\lambda $
determines\footnote{In fact, we must choose a particular representative of
the cohomology class~$\lambda \in H^n(BG;\ZZ)$ to carry out the constructions
which follow~\cite{FQ, Appendix~B}, ~\cite{F1,\S2}.} a characteristic
number~$\lambda _X(P)$ for a $G$~bundle $P\to X$ over a closed oriented
$n$-manifold~$X$.  One can define a $\ZZ$-torsor $\tors YQ$ for a $G$~bundle
$Q\to Y$ over a closed oriented $(n-1)$-manifold~$Y$, and a relative
invariant for $n$-manifolds with boundary.  The story continues to higher
codimensions.  For example, the invariant of a bundle over a closed oriented
$(n-2)$-manifold is a certain type of category, called a {\it gerbe\/}.

We illustrate the codimension~2 invariant in the simplest case of the Euler
number.  Recall the definition~\thetag{1.3} of the $\ZZ$-torsor associated to
a line bundle $K\to Y$ over a closed oriented surface.  Now suppose $Y$~is a
compact oriented 1-manifold with boundary, i.e., a finite union of circles
and closed intervals.  Let $u\:\bY\to \bK$ be a nonzero section over the
boundary and let $\Sect^*(K,u)$~be the set of nonzero sections $t\:Y\to K$
with~$\partial t=u$. Then define the $\ZZ$-torsor
  $$ T_Y(K,u) = \{e\:\Sect^*(K,u)\to Z : e(t')-e(t)= \deg_{\bX}(t'/t)\}.
      $$
The degree is well-defined since~$\partial t'=\partial t$.  We need to
determine the dependence of~$T_Y(K,u)$ on~$u$.  Now if $u'$~is any other
trivialization of $\bK\to \bY$, then the set of nonzero paths joining~$u$
to~$u'$, up to homotopy, is a $\ZZ$-torsor $T_{\bY}(\bK,u,u')$.  Such a path
gives a 1:1~correspondence $\Sect^*(K,u)\to \Sect^*(K,u')$, and so there is
an isomorphism
  $$ T_Y(K,u)\times T_{\bY}(\bK,u,u')\longrightarrow T_Y(K,u').  $$
This leads us to define
  $$ \multline
      \sG_{\bY}(\bK) = \bigl\{T\:\Sect^*(\bK)\to \text{$\{\ZZ\text{-torsors}\}$
     with given isomorphisms}\\
      \text{$T(u)\times T_{\bY}(\bK,u,u') \to T(u')$ which satisfy various
     conditions} \bigr\}. \endmultline\tag{1.9} $$
In this definition $T$~is a set-valued function, so $\sG_{\bY}(\bK)$~is a
collection of set-valued functions, in particular, a category.  The ``various
conditions'' are related to the category structure and describe an
associativity constraint for the isomorphism associated to three nonzero
sections~$u,u',u''$.  Think of the {\it gerbe\/}~\thetag{1.9} as an affine
copy of the category of $\ZZ$-torsors.  The relative invariant~$T_Y(K)$ lies
in the gerbe~$\sG_{\bY}(\bK)$ by definition.

This invariant admittedly contains little information.  However, if we start
with a 4~dimensional characteristic class, then the gerbe attached to a
surface leads to a central extension of the diffeomorphism group of the
surface which arises in quantum Chern-Simons theory~\cite{W1}.  (It is better
here to start with the {\it signature\/} of a 4-manifold~\cite{F2}.)  The
central extension is often realized in terms of ``2-framings''.

Our interest in~\S{2} is in 3~dimensional invariants of principal bundles
with {\it finite\/} gauge group~$G$.  But if $G$~is finite then $H^n(BG;\ZZ)$
consists entirely of torsion elements for~$n>0$, and so the primary integral
characteristic numbers vanish.  Rather, we start with
  $$ \lambda \in H^3(BG;\RR/\ZZ) \cong  H^4(BG;\ZZ) \tag{1.10} $$
and obtain a characteristic number $\lambda _X(P)\in \RZ$ for $G$~bundles
$P\to X^3$ over a closed oriented 3-manifold.  Exponentiating, the invariant
  $$ \el P\in \TT \tag{1.11} $$
lies in the circle group~$\TT$ of unit norm complex numbers.  The relative
picture is similar: To a $G$~bundle $Q\to Y$ over a closed oriented
surface~$Y$ we attach a $\TT$-torsor, and to a bundle $P\to X$ over a
3-manifold with boundary the relative invariant lives in the $\TT$-torsor of
the boundary.  These invariants obey a gluing law~\thetag{1.6} which we now
write multiplicatively.  Notice that any $\TT$-torsor is the set of unit norm
elements in a hermitian line, which here we denote~$L_Y(Q)$.  So the relative
invariant of~$P\to X$ is
  $$ \el P\in L_{\bX}(\bP). \tag{1.12} $$
Similarly, the invariant of a $G$~bundle $R\to S$ over a closed oriented
1-manifold~$S$ (finite union of circles) is a $\TT$-gerbe~$\sG(R)$, which we
think of as an affine copy of the category of hermitian lines.  (Imagine
definition~\thetag{1.9} with hermitian line-valued functions.)  Then the
relative invariant of a $G$~bundle $Q\to Y$ over an oriented surface with
boundary is an element
  $$ L_Y(Q)\in \sG_{\bY}(\bQ). \tag{1.13} $$

Although~\thetag{1.11} is a topological invariant, and in that sense is
primary (Table~1), because of the transgression in~\thetag{1.10} we can think
of it as a Chern-Simons invariant for bundles with finite structure group.
More generally, for any compact group~$G$ and class~$\lambda \in H^n(BG;\ZZ)$
there are secondary {\it geometric\/} invariants for connections on bundles
over compact oriented manifolds of dimension at most~$n-1$.  For a closed
oriented $(n-1)$-manifold the invariant lies in~$\TT$, for a closed oriented
$(n-2)$-manifold the invariant {\it is\/} a $\TT$-torsor (or hermitian line),
etc.  For a family of connections on an $(n-2)$-manifold there is a
connection on the line bundle over the parameter space.  We develop these
ideas in~\cite{F3} using a variant of the {\it \v Cech-de Rham complex\/}
(more fondly ``tic-tac-toe''~\cite{BT}) which is related to Deligne
cohomology~\cite{Bry}.  This amounts to a generalization of the theory of
differential forms, including both differentiation and integration.  The case
$n-1=3$ is of special interest~\cite{F4} because of the relation with special
geometric structures in low dimensional gauge theory and with the
corresponding quantum invariants.

Let us briefly consider the $K$-theory side of Table~1.  From the point of
view of Dirac operators the primary characteristic number---in~$\ZZ$---is the
index of a Dirac operator on an even dimensional closed spin manifold, and
the secondary characteristic number---in~$\TT$---is the exponentiated $\xi
$-invariant of a Dirac operator on an odd dimensional closed spin manifold.
(Recall that the $\xi $-invariant is half the $\eta $-invariant plus half the
dimension of the kernel.)  Notice that the former is a topological invariant
whereas the latter depends on geometric data, i.e., a metric and possibly a
connection.  One can define relative invariants in this context.  First, a
reinterpretation of the work of Atiyah-Patodi-Singer~\cite{APS} constructs
the index of a Dirac operator on a manifold with boundary as a topological
invariant~\cite{F5}.  It lives in a $\ZZ$-torsor (which is an invariant of
the Dirac operator on the boundary) and satisfies a gluing law.  This is
completely analogous to~\thetag{1.4} and~\thetag{1.6}.  In recent work with
Xianzhe Dai~\cite{DF} we construct a relative exponentiated $\xi $-invariant
(analogous to~\thetag{1.12}) which lives in the determinant line of the Dirac
operator on the boundary and satisfies a gluing law.  However, I do not know
how to extend these constructions to manifolds of lower dimension (except for
the next step on the ``primary side'').  Also, it would be interesting to
give topological constructions of the relative primary invariants which
generalize the direct image~\cite{AH}.

\newpage
\head
\S{2} Generalized Path Integrals
\endhead
\comment
lasteqno 2@ 10
\endcomment

Consider a finite group~$G$ and a class $\lambda \in H^3(BG;\RZ)$.  In the
last section we indicated that there are invariants~$\lambda _X(P)\in \TT$
for each principal $G$~bundle $P\to X^3$ over a closed oriented 3-manifold.
If we are interested in constructing topological invariants of~$X$, then we
need to eliminate the dependence on~$P$.  One way---the physicists' way---is
to integrate over~$P$.  In quantum field theory $P$~is called a {\it field\/}
and this integral over the space of fields is called the {\it path
integral\/}.  It is in this sense that our discussion here is a field theory.
More specifically, it is the quantum Chern-Simons theory~\cite{W1} for finite
gauge group, first introduced by Dijkgraaf and Witten~\cite{DW}.
Fortunately, the integral reduces to a finite sum in this theory and there
are no analytic problems to worry about.

Let $\fld X$~denote the collection of $G$~bundles over~$X$.  It is a
category---morphisms are bundle maps which cover the identity.  Let $\fldb
X$~denote the finite set of equivalence classes of bundles, and define a
measure on~$\fldb X$ by
  $$ \mu _X(P) = \frac{1}{\#\Aut P} \tag{2.1} $$
for $P\in \fld X$.  The topological invariant of a closed oriented
3-manifold~$X$ we wish to study is
  $$ Z_X = \int_{\fldb X} \el P \,d\mu _X(P). \tag{2.2} $$
We view the integrand~\thetag{1.11} as a complex number and perform the
integral by summing complex numbers.  So $Z_X\in \CC$.  Although we use
integral notation, \thetag{2.2}~is a finite sum.  The fact that $\el P$~is a
topological invariant of~$P$ quickly leads to a proof that $Z_X$~is a
topological invariant of~$X$.


Just as with the Euler number in~\S{1}, we would like to define an invariant
of a 3-manifold with boundary which satisfies a gluing law analogous
to~\thetag{1.6}.  Let $X$~be a closed oriented 3-manifold and
$Y\hookrightarrow X$ an embedded closed oriented surface.  Denote by~$X\cut$
the manifold obtained by cutting~$X$ along~$Y$.  (See Figure~1.)  Notice that
$\partial X\cut = Y \sqcup -Y$ is a disjoint union of two copies of~$Y$.  The
fields (bundles) fit into the following diagram:
  $$ \CD
     \fld X @>c>> \fld{X\cut}\\
     @Vr_1VV @VVr_2V\\
     \fld Y @>\Delta >> \fld Y\times \fld{-Y}
     \endCD  $$
The vertical arrow~$r_1$ is restriction to~$Y$, the arrow~$r_2$ is
restriction to~$\partial X\cut$, the arrow~$\Delta $ is the diagonal
inclusion, and $c$~is the pullback under the gluing map.  For the moment we
ignore symmetries and pretend that \thetag{2.2}~is an integral over~$\fld X$.
Then we propose to do the integral over~$\fld X$ in two stages using Fubini's
theorem: First integrate over the fibers of~$r_1$ and then over~$\fld Y$.
Now the gluing law~\thetag{1.9} for~$\lambda $ says that if $P\in \fld X$ and
$P\cut=c(P)\in \fld{X\cut}$, then
  $$ \el P = e^{2\pi i\lambda _{X\cut}(P\cut)}. \tag{2.3} $$
The right hand side of~\thetag{2.3} lives in $L_Y(Q)\otimes L_{-Y}(Q)$
(cf.~\thetag{1.12}) which is identified with~$\CC$ via a pairing analogous
to~\thetag{1.5}.  Our hope, then, is to make the following computation:
  $$ \split
      Z_X=\int_{\fld X}\el P \,d\mu _X(P) &= \int_{\fld Y}\int_{r_1\inv (Q)}
     \el P \,d\mu _{r_1\inv (Q)}(P) \, d\mu _Y(Q)\\
      &= \int_{\fld Y} \int_{r_2\inv (Q,Q)} e^{2\pi i\lambda _{X\cut}(P\cut)}
     \,d\mu _{r_2\inv (Q,Q)}(P\cut) \, d\mu _Y(Q).\endsplit \tag{2.4} $$
For this to be a valid computation we need the measures to work out properly.
Also, we must include the symmetries.  Both are easily handled~\cite{FQ,\S2}.

Let's reinterpret the last line of~\thetag{2.4}.  For a moment suppose that
$X'$~is any compact oriented 3-manifold with boundary.  Then for $Q\in
\fld{\bX'}$ define
  $$ \fld{X'}(Q) = \{P\to X' \text{ such that $\partial P=Q$} \}.
     $$
There are symmetries here as well, and we are deliberately vague in order to
keep the ideas as simple as possible.  Note that $\fld{X'}(Q)$ is the space
of fields with a given fixed boundary value.  We generalize the path
integral~\thetag{2.2} to manifolds with boundary by defining
  $$ Z_{X'}(Q) = \int_{\fldb{X'}(Q)} e^{2\pi i\lambda _{X'}(P)}\,d\mu
     _{X'}(P).  \tag{2.5} $$
The right hand side takes values in the hermitian line~$L_{\bX'}(Q)$
(cf.~\thetag{1.12}).  So $Z_{X'}$~is a section of the hermitian line bundle
$L_{\bX'}\to \fld{\bX'}$.  In an appropriate sense it is invariant under
symmetries.  For any closed surface~$Y$ set
  $$ E(Y) = \text{invariant sections of $L_{Y} \to \fld{Y}$.}
     \tag{2.6} $$
Then \thetag{2.5} determines a relative invariant
  $$ Z_{X'} \in E(\bX').  $$
We impose an $L^2$~inner product on~\thetag{2.6} using the
measure~\thetag{2.1} on {\it equivalence classes\/} of bundles.

With these definitions the gluing law~\thetag{2.4} (with symmetries restored)
takes the form
  $$ \split
      Z_X &= \int_{\fldb {\bX}} Z_{X\cut}(Q,Q) \,d\mu _{\bX}(Q)\\
      &= \Trace_{\bX}(Z_{X\cut}),\endsplit \tag{2.7} $$
where $\Trace_{\bX}\:E(\bX)\otimes E(-\bX)\to\CC$ is formed using the
$L^2$~inner product.  (Notice that $L_{-\bX}(Q)\cong \overline{L_\bX(Q)}$ and
so $E(-\bX)\cong \overline{E(\bX)}$.)

What we have recounted so far is the standard argument that path integrals
give numerical invariants which satisfy gluing laws.  The relative invariants
in quantum theories do not live in ``one dimensional'' torsors, as in the
classical case~\thetag{1.3}, \thetag{1.8}, but rather in Hilbert
spaces~\thetag{2.6}.  Now we want to go further and derive gluing laws for
the Hilbert spaces.

Here is the main idea: Re-express the Hilbert space~$E(Y)$ as an integral
over the space of fields~$\fldb Y$.  Then repeat the argument which leads
from~\thetag{2.2} to the gluing law~\thetag{2.7}.  In the process we will
define an invariant~$\sE(S)$ of a closed oriented 1-manifold~$S$, by analogy
with~\thetag{2.6}.  If the latter can be re-expressed as an integral we can
again iterate the process.

To find the integral which computes~$E(Y)$ recall that the characteristic
class~$\lambda $ associates to each $G$~bundle $Q\to Y$ a hermitian
line~$L_Y(Q)$.  Furthermore, an isomorphism $Q'\to Q$ of bundles induces an
isomorphism $L_Y(Q')\to L_Y(Q)$ of hermitian lines.  So there is a quotient
hermitian line bundle
  $$ L_Y\longrightarrow \fldb Y  $$
(This bundle degenerates where automorphisms of~$Q\in \fld Y$ act
nontrivially on~$L_Y(Q)$.)  Our formula for the Hilbert space is:
  $$ \boxed{E(Y) = \int_{\fldb Y} L_Y(Q)\,d\mu _Y(Q)} \tag{2.8} $$
This is a boldly presented formula and it requires some explanation.  The
right hand side has the form
  $$ \sum\limits_{\text{finite} } (\text{positive number} )\cdot
     (\text{hermitian line} ). \tag{2.9} $$
If $\mu >0$ and $L$~is a hermitian line, define $\mu \cdot L$ to be the
hermitian line with the same underlying vector space as~$L$ but with an inner
product which is $\mu $~times the inner product on~$L$.  Interpret the finite
sum as the direct sum of hermitian lines.  This is the sense in which the
right hand side of~\thetag{2.8} defines a Hilbert space---as a sum of
hermitian lines---and it is easy to see that \thetag{2.6}~and
\thetag{2.8}~define the same Hilbert space.  As explained in~\S{1} the
hermitian line~$L_Y(Q)$ is the two dimensional counterpart of the three
dimensional invariant~$\el P$.  So the integrals~\thetag{2.2}
and~\thetag{2.8} have the same form: they are the integral of the
exponentiated classical action over the space of fields.

By now it should be clear how to continue.  If $Y$~is an oriented surface
with boundary, then by analogy with~\thetag{2.5} we define
  $$ E(Y)(R) = \int_{\fldb Y(R)} L_Y(Q)\, d\mu _Y(Q),\qquad R\in \fld{\bY}.
      $$
Recall from~\thetag{1.13} that $L_Y(Q)$~is an element of the
gerbe~$\sG_{\bY}(R)$.  Think of the latter as an affine copy of the category
of hermitian lines.  Then just as we can add hermitian lines to obtain a
Hilbert space~\thetag{2.9}, we can add elements of~$\sG_{\bY}(R)$.  What do
we get?  If $\sG_{\bY}(R)$~were trivial we would get a Hilbert space.  Then
$E(Y)$~would be an invariant function from~$\fld{\bY}$ to the category of
Hilbert spaces.  The nontriviality of~$\sG_{\bY}(R)$ means that $E(Y)$~is an
invariant function from~$\fld{\bY}$ to an affine version of the category of
Hilbert spaces, a so-called {\it 2-Hilbert space\/}~\cite{KV}.  The
space~$\sE(\bY)$ of all such functions is again a 2-Hilbert space and
$E(Y)\in \sE(\bY)$.  If we cut a surface along a circle, then we can
formulate a gluing law analogous to~\thetag{2.7}.  Finally, for any closed
oriented 1-manifold~$S$ we write the 2-Hilbert space~$\sE(S)$ as an integral:
  $$ \sE(S) = \int_{\fldb S}\sG_S(R)\,d\mu _S(R). \tag{2.10} $$
In principle, we can continue to even lower dimensions.

We stop here, but refer the reader to~\cite{F2} for an expository account
which begins with the notion of a 2-Hilbert space.  There we give a general
argument to show how in a 3~dimensional topological quantum field theory a
``quantum group'' arises from the 2-Hilbert space.  In~\cite{F1} we carry out
detailed computations for gauge theory with finite gauge group.

We end with some remarks about Chern-Simons theory with {\it continuous\/}
gauge group.  Then \thetag{2.2}, \thetag{2.8}, and~\thetag{2.10} are
replaced by integrals over the space of connections mod equivalence.
Equation~\thetag{2.2} is then the path integral heuristic given by
Witten~\cite{W1}.  Reshetikhin and Turaev~\cite{RT} subsequently gave an
explicit computable formula in terms of quantum group data and proved that it
is a topological invariant.  Note that the formal integral~\thetag{2.2} is
over an infinite dimensional space.  By contrast, the integral~\thetag{2.8}
reduces, after symplectic reduction, to an integral over the finite
dimensional moduli space~$\sM_Y$ of flat connections on~$Y$.  The integrand
is a line bundle with connection whose curvature is the symplectic form
on~$\sM_Y$.  In this situation we may imagine that \thetag{2.8}~is a formal
expression for geometric quantization.  The use of an integral is justified
by the existence of gluing laws---Verlinde's formula~\cite{V}.  By now there
are various proofs of at least special cases of these gluing laws.

What about~\thetag{2.10}?  If $S$~is a circle then up to equivalence a
connection is determined by its holonomy, which is a conjugacy class in the
gauge group~$G$.  So \thetag{2.10}~is an integral over a compact Lie
group~$G$ with some invariance under the adjoint action.  The integrand is a
geometric object over~$G$, ``invariant'' under the adjoint action, with a
``connection'' whose curvature is the canonical biinvariant 3-form on~$G$
constructed from the starting data $\lambda \in H^4(BG)$.  (See~\cite{Bry}
for one description.)  I think it is reasonable to imagine that there is a
geometric process---analogous to geometric quantization---which constructs
the 2-Hilbert space~$\sE(S)$ from this ``gerbe bundle with connection'', but
for now this remains a mystery.

\enddocument